\documentclass[12pt,a4paper]{amsart}
\usepackage{amssymb}
\makeatletter
\renewcommand{\@begintheorem}[2]{
\rm \trivlist \item [\hskip \labelsep {\bf #2\ \ #1.}]
                                }
\makeatother

\makeatletter

\newcommand{\ZZ}{{\bf Z}}
\newcommand{\QQ}{{\bf Q}}
\newcommand{\RR}{{\bf R}}
\newcommand{\CC}{{\bf C}}
\newcommand{\DD}{{\bf D}}

\newcommand{\HH}{{\mathbb H}}
\newcommand{\PP}{{\bf P}}

\newcommand{\BB}{{\mathbb B}}

\newcommand{\cA}{{\mathcal A}}

\newcommand{\cF}{{\mathcal F}}

\newcommand{\cT}{{\mathcal T}}
\newcommand{\cD}{{\mathcal D}}
\newcommand{\cP}{{\mathcal P}}

\newcommand{\cX}{{\mathcal X}}

\newcommand{\cO}{{\mathcal O}}

\newcommand{\ra}{\rightarrow}

\newcommand{\rmd}{\mbox{d}}

\newcommand{\bes}{\begin{equation*}}
\newcommand{\ees}{\end{equation*}}

\title{Examples of Calabi-Yau threefolds parametrised by Shimura varieties}
\author{Alice Garbagnati}
\author{Bert van Geemen}
\address{Dipartimento di Matematica, Universit\`a di Milano,
Via Saldini 50, I-20133 Milano, Italia}
\email{alice.garbagnati@unimi.it}
\email{lambertus.vangeemen@unimi.it}

\begin{document}

\begin{abstract}
These are notes from talks of the authors on some explicit examples of families of Calabi-Yau threefolds which are parametrised by a Shimura variety. We briefly review the periods of Calabi-Yau threefolds and we discuss a recent result on Picard-Fuchs equations for threefolds which are hypersurfaces with many automorphisms. Next various examples of families parametrised by Shimura varieties are given. Most of these are due to J.C.\ Rohde. The examples with an automorphism of order three are given in some detail. We recall that such families do not have maximally unipotent monodromy and that the Shimura varieties in these cases are ball quotients.
\end{abstract}

\maketitle

Calabi-Yau threefolds have been intensively studied in the context of 
Mirror Symmetry, which originated in theoretical physics. 
The most basic version of Mirror Symmetry states that given a CY threefold
$X$, there should exist a CY threefold $Y$, the Mirror of $X$, 
with $h^{1,1}(X)=h^{2,1}(Y)$ and
$h^{2,1}(X)=h^{1,1}(Y)$. A further requirement is that
the variation of Hodge structures on the third cohomology group in the deformations 
of $X$ is related, in specific way, to the K\"ahler cone in the second cohomology group of 
the deformations of $Y$. There are now quite a few cases where such Mirror pairs have 
been found and where profound aspects of Mirror Symmetry, like the relation with Gromow-Witten invariants, have been verified. 

As CY threefolds are assumed to be projective (or at least K\"ahler), in case
$h^{2,1}(X)=0$ there cannot exist a Mirror of $X$. Quite a few of such (rigid) threefolds $X$ are known and a modification of Mirror Symmetry, in one specific case, was proposed in \cite{CDP}. 
The expected relation between deformations of $X$ and the K\"ahler cone of $Y$ requires
that there exist boundary points in the (complex structure) moduli space of $X$ where the variation of Hodge structures on the $H^3$ has maximally unipotent monodromy. 
It was recently pointed out by J.C.\ Rohde 
\cite{Rohde-maximal} that there do exist families of CY threefolds which do not admit 
such boundary points. The moduli spaces of the families in question are Shimura varieties.
No modification of Mirror Symmetry in these cases is known to us.

Having a base which is a Shimura variety is otherwise a quite desirable property because then one has a very good control over the variation of the Hodge structures. Moreover, the so-called CM points (cf.\ \cite{Bcycm}, \cite{Rlnm}) will be dense in the moduli space. Physicists expect the field theory
on the corresponding CY threefolds to be simpler \cite{GV}.

In these notes, we give various examples of CY threefolds whose moduli space is a Shimura variety. A very simple example, where the Shimura variety is the moduli space of elliptic curves, is given in section \ref{easy}. The CY threefolds are of Borcea-Voisin type
and it is easy to see that the family does not have maximally unipotent monodromy. 
In section \ref{e3} we discuss various examples, due to J.C.\ Rohde, of CY threefolds with an automorphism of order three.
Their moduli spaces are Shimura varieties which are ball quotients and they do not have maximally unipotent monodromy. 
We show, in all cases but one, that Rohde's CY threefolds are desingularisations of a quotient of a product of two fixed elliptic curves with another curve of higher genus. These curves of higher genus  also have an automorphism of order three. The moduli space of CY threefolds is, locally at least,
the moduli space of such pairs $(C,\beta)$ where $C$ is the higher genus curve and 
$\beta$ is its automorphism of order three. This allows one to write down the Picard-Fuchs equation for the variation of Hodge structures explicitly in one case, see section \ref{exR}.
New examples similar to these, but using an automorphism of order four, can be found in \cite{Ga}.

\section{Periods of Calabi-Yau threefolds.}

\subsection{} Good references for CY threefolds and Mirror Symmetry are the overview of M.\ Gross in \cite{GHJ} and the book \cite{CK}.

\subsection{Calabi-Yau threefolds}\label{polhs}
In these notes, a Calabi-Yau threefold $X$ is a smooth, three dimensional 
(complex) projective variety with trivial canonical bundle,
$$
\Omega^3_X\,\cong\,\cO_X,\qquad\mbox{and}\qquad H^1(X,\cO_X)\,=\,0,
$$
using Serre duality one then finds 
$H^2(\cO_X)\cong H^1(\Omega_X^3)^*=H^1(\cO_X)=0$ and 
$H^3(\cO_X)\cong H^0(\Omega_X^3)^*\cong \CC$.

Examples of CY threefolds are quintic hypersurfaces in $\PP^4$ and complete
intersections of two hypersurfaces of degree three in $\PP^5$. 
Many more families of CY threefolds can be found as hypersurfaces in four dimensional toric varieties \cite{KS}.
There are obvious restrictions on the topology of $X$, in fact Hodge theory easily implies that 
$H^1(X,\CC)=0$ and $\dim H^3(X,\CC)\geq 2$.
As $X$ is projective, one has $\dim H^2(X,\CC)\geq 1$. 
Poincar\'e duality implies that $b_i=b_{6-i}$, but
no further restrictions on the Betti numbers $b_i:=\dim H^i(X,\CC)$ of $X$ are known.

In analogy with the case of curves, abelian varieties and K3 surfaces, 
one studies the Hodge structures on the cohomology groups
in order to understand these varieties better. For CY threefolds, only $H^3$ is of interest, as $H^2(X,\CC)=H^{1,1}(X)$. 

\subsection{The polarised Hodge structure on $H^3$}
The Hodge structure on $H^3(X,\ZZ)/\mbox{torsion}$
is the decomposition of its complexification:
$$
H^3(X,\CC)\,=\,
\underbrace{
\underbrace{H^{3,0}(X)}_{F^3}\,\oplus\,H^{2,1}(X)}_{F^2}
\,\oplus\,H^{1,2}(X)\,\oplus
H^{0,3}(X),\qquad H^{p,q}(X)\,=\,\overline{H^{q,p}}(X).
$$
From the $F^2$ of the Hodge filtration, one recovers
$
H^3(X,\CC)=F^2\,\oplus\,\overline{F^2}
$.
The intersection form on $H^3(X,\ZZ)$, which factors over 
$H^3(X,\ZZ)/\mbox{torsion}$, defines a polarisation $Q_X$ on this Hodge structure:
$$
Q_X:\,(H^3(X,\ZZ)/\mbox{torsion})\,\times\,(H^3(X,\ZZ)/\mbox{torsion})\,
\longrightarrow\,\ZZ,\qquad
Q_X(\theta_1,\theta_2)\,:=\,i^3\int_X\theta_1\wedge\overline{\theta_2}
$$
where $i=\sqrt{-1}\in\CC$ and we identified $H^3(X,\ZZ)/\mbox{torsion}$ with the
image of $H^3(X,\ZZ)$ in $H^3(X,\RR)\cong H^3_{DR}(X)$, the de Rham cohomology group.
The polarisation $Q_X$ is a symplectic form (so it is non-degenerate, unimodular
and alternating). 
It extends to a Hermitian form $H_X:=iQ_X$ on $H^3(X,\CC)$ 
for which the Hodge decomposition is orthogonal:
$$
H_X(v,w)\,=\,0\qquad\mbox{if}\quad v\in H^{p,q}(X),\;w\in H^{r,s}(X)
\quad\mbox{and}\quad (p,q)\neq (r,s)
$$
and which is positive/negative definite on the $H^{p,q}(X)$:
$$
H_X\,:=\,iQ_X\quad\mbox{is}\quad \left\{\begin{array}{rrrr}
>0&\mbox{on}\quad H^{3,0}(X),&(<0&\mbox{on}\quad H^{0,3}(X)),\\
<0&\mbox{on}\quad H^{2,1}(X),&(>0&\mbox{on}\quad H^{1,2}(X)).
\end{array}\right.
$$

\subsection{The Period domain}\label{pdom}
Let $N=b_3$ be the rank of $H^3(X,\ZZ)$ and let $Q$ be 
a (fixed) symplectic form on $V_\ZZ:=\ZZ^N$. Then we consider the period space $\cD=\cD_N$ of all polarized weight three Hodge structures on $(V_\ZZ,Q)$
of CY-type. 
An element of $\cD$ is a decomposition 
$$
V_\CC\,:=\,V_\ZZ\otimes_\ZZ\CC\,=\,V^{3,0}\oplus V^{2,1}\oplus V^{1,2}\oplus V^{0,3},\qquad V^{p,q}\,=\,\overline{V^{q,p}}
$$
such that the Hermitian form $H(v,w):=iQ(v,\overline{w})$, with $v,w\in V_\CC$
and where $Q$ is extended $\CC$-linearly, 
is positive definite on $V^{3,0},V^{1,2}$ and negative definite on $V^{2,1},V^{0,3}$. Moreover we require $\dim V^{3,0}=1$ and we denote $q:=\dim V^{2,1}$ (so $N=2+2q$).

The $q+1$-dimensional subspaces $W$ of $V_\CC$ on which $H$ is positive
are parametrized by a $(q+2)(q+1)/2$ dimensional variety (it is isomorphic to the Siegel half space of $(q+1)\times(q+1)$ complex symmetric matrices with positive definite imaginary part). The one dimensional subspaces $W_1$ of such a subspace $W$ are the points of $\PP W$ and thus are parametrized by a $q$-dimensional variety. Given $W_1\subset W$, let $W_1^\perp$ be the orthogonal complement of $W_1$ in $W$. Then we obtain a polarized Hodge structure on $V_\ZZ$ by defining: 
$$
V^{3,0}\,:=\,W_1,\qquad V^{2,1}\,:=\,\overline{W_1^\perp},
\qquad V^{1,2}\,:=\,W_1^\perp,\qquad V^{0,3}\,:=\,\overline{W_1}.
$$
Conversely, any polarized Hodge structure on $(V_\ZZ,Q)$ defines a complex line $W_1$ in a positive $(q+1)$-dimensional subspace $W$, hence we get 
$$
\dim \cD\,=\,q\,+\,(q+2)(q+1)/2\,=\,(q^2+5q+2)/2,\qquad 
q:=\dim V^{2,1}.
$$

\subsection{The Period map}
A marking of the CY threefold $X$ is a symplectic
isomorphism $(H^3(X,\ZZ)/\mbox{torsion})\rightarrow V_\ZZ$. 
The 
$\CC$-linear extension of this isomorphism maps the
Hodge decomposition of $H^3(X,\CC)$
to a decomposition of $V_\CC$. In this way we obtain a polarized Hodge structure on $V_\ZZ$. 
In particular, we get a point $\cP(X)\in\cD$, the period point of $X$.

\subsection{Deformations of $X$}\label{defs}
An important result, due to Bogomolov, Tian and Todorov, 
on CY varieties is that the deformations are unobstructed.
The first order deformations of a complex variety are parametrised by
the cohomology group $H^1(X,\cT_X)$, where $\cT_X$ is the holomorphic tangent bundle of $X$. As $X$ is a CY threefold, the cup product pairing 
$\Omega^1_X\times\Omega^2_X\rightarrow \Omega^3_X\cong\cO_X$
gives a duality $\cT_X\cong (\Omega^1_X)^*\cong \Omega^2_X$
and thus $H^1(X,\cT_X)\cong H^1(X,\Omega^2_X)\cong H^{2,1}(X)$.
The unobstructedness asserts that there is a 
neighbourhood $B$ of $0\in H^1(X,\cT_X)$ and there is a 
family of CY threefolds $\pi:\cX\rightarrow B$ with fiber $\pi^{-1}(0)=X$, such that the  period map $\cP:B\rightarrow \cD$ has an injective 
differential:
$$
(\rmd \cP)_0\,:\,T_0B\,=\,H^1(X,\cT_X)\,\cong\,H^{2,1}(X)\,
\longrightarrow\, T_{\cP(X)}\cD\qquad\mbox{is injective}.
$$
Here we used Ehresmann's theorem which asserts that,
if $B$ is chosen small enough,
there is a diffeomorphism $\phi:\cX\rightarrow B\times X$ such that $\pi_B\phi=\pi$.
As $\phi$ induces isomorphisms 
$H^3(X_b,\ZZ)\cong H^3(X,\ZZ)$ for any $X_b:=\pi^{-1}(b)$,
we can extend the marking on $X$ to a marking 
$(H^3(X_b,\ZZ)/\mbox{torsion})\rightarrow V_\ZZ$
and thus we get the period map $\cP:B\rightarrow \cD$.

Any family of CY threefolds which contains $X$ is locally near $X$ obtained as the pull-back from the family $\cX\rightarrow B$. 
Therefore the image of the period map of any family has dimension at most
$q=\dim H^{2,1}(X)$ and the image of $B$ has codimension
$(q+2)(q+1)/2$ in $\cD$. Recent studies of the geometry of $\cD$ and  these subvarieties are \cite{CGG} and \cite{LSY}.

\subsection{The Picard-Fuchs equation}
In this section we will assume for simplicity that $q=\dim H^{2,1}(X)=1$.
With the notation of section \ref{defs}, 
the diffeomorphism $\phi:\cX\rightarrow B\times X$ induces an isomorphism
of sheaves 
$R^3\pi_*\ZZ\stackrel{\cong}{\rightarrow} H^3(X,\ZZ)_B$ on $B$, 
where the last sheaf is just the locally constant sheaf defined by the 
abelian group $H^3(X,\ZZ)$. 
Using the marking $H^3(X,\ZZ)/\mbox{torsion}\cong V_\ZZ$ 
and the Hodge decomposition of $H^3(X_t,\CC)$ for each deformation $X_t$ of $X$, we obtain a (trivial) vector bundle $V_\CC\times B$ over $B$
with holomorphic subbundles 
$$
\cF^3\,\subset\,\cF^2\,\subset\,\cF^1\,\subset\cF^0
\cong V_B:=V_\CC\times B,
\qquad \cF^3_t=H^{3,0}(X_t),
$$
where we identify $V_\CC\times\{t\}$ with $H^3(X_t,\CC)$.

The period map $\cP$ describes the variation of these subbundles inside the trivial bundle $V_\CC\times B$. Another way to describe this variation is to
take a non-vanishing section $\omega$ of the rank one bundle $\cF^3$,
so $\omega(t)$ is a basis of $H^{3,0}(X_t)$ for all $t\in B$.
The trivial bundle $V_B$ comes with the Gauss-Manin connection $\nabla$
which maps the horizontal sections 
$s_v:=t\mapsto (v,t)$ to zero, where $v\in V_\CC$:
$$
\nabla\,=\,\nabla_{\partial/\partial t}:\;
V_\CC\times B\,\longrightarrow\, V_\CC\times B.
$$
Applying the connection 
$i$ times to the section $\omega$, we get a section $\nabla^i\omega$.
As $\dim V_\CC=2+2q=4$, there must be a linear relation, with coefficients
$p_i(t)$ which will be holomorphic in $t$:
$$
\DD\,\omega\,=\,0,\qquad \DD\,:=\,\sum_{i=0}^4p_i(t)\nabla^i.
$$
This linear relation is known as the Picard-Fuchs equation.

Instead of considering this rank four bundle with its section $\omega$, one can also choose a basis $\gamma_1,\ldots,\gamma_4$ of 
$H_3(X,\ZZ)/\mbox{torsion}$ and define four holomorphic functions
$\varphi_i(t):=\int_{\gamma_i}\omega(t)$ on $B$,
where $\gamma_i$ is identified with a cycle in $H_3(X_t,\ZZ)/\mbox{torsion}$ using the diffeomorphism $\phi$.
These four functions are a basis of the solutions of 
the degree four differential  operator 
$\sum_{i=0}^4p_i(t)(\rmd/\rmd t)^i$
which is also called the Picard-Fuchs equation for the family 
$\cX\rightarrow B$.

\subsection{An example}
The Dwork pencil of quintic threefolds in $\PP^4$
is defined by the equation
$$
X_t:\qquad X_1^5+\ldots+X_5^5\,-\,5tX_1X_2\cdots X_5\,=\,0.
$$
For general $t\in\CC$, the variety $X_t$ is a CY threefold with 
$h^{1,1}(X_t)=1$ and $q=h^{2,1}(X_t)=101$. However, there is a finite subgroup $G\cong (\ZZ/5\ZZ)^3$ acting on $\PP^4$ which induces automorphisms on each $X_t$ and the third cohomology group splits under this action:
$$
H^3(X_t,\QQ)\,=\,T_t\,\oplus\,S_t,\qquad 
T_t\,:=\,H^3(X_t,\QQ)^G\,\cong\,\QQ^4,\qquad H^{3,0}(X_t)\subset T_t\otimes_\QQ\CC.
$$
Thus the $G$-invariant part of the cohomology gives a four dimensional
variation of polarised Hodge structures and thus it gives a degree four Picard-Fuchs equation.

In the context of Mirror Symmetry, it was observed that the (singular)
quotient variety $X_t/G$ has a resolution of singularities $M_t$ 
which is a CY threefold, moreover its Hodge numbers are:
$$
h^{1,1}(M_t)\,=\,101,\qquad h^{2,1}(M_t)\,=\,1,\qquad
M_t\,:=\,\widetilde{X_t/G}.
$$
Note that $h^{p,q}(M_t)=h^{3-p,q}(X_t)$, which is one of the requirements for the ($101$-dimensional) family of quintic CY threefolds and the one parameter family of $M_t$'s to be Mirror CY families.

The quotient map induces an isomorphism 
$T_t\cong H^3(M_t,\QQ)$. In particular, the degree four Picard 
Fuchs equation obtained from the variation of the $T_t$ is the Picard-Fuchs 
equation of the one parameter family of CY threefolds $M_t$. 
A spectacular result from Mirror Symmetry is that a certain solution of this Picard-Fuchs equation 
defines a power series in one variable whose coefficients $a_d$ allow one to
compute the Gromow-Witten invariants of a quintic threefold, that is, roughly, the  number of rational curves of degree $d$ on a quintic threefold.

In the paper \cite{GPR}, Greene, Plesser and Roan verify that 
there is an action of the group $H\cong\ZZ/41\ZZ$ on $\PP^4$ such that
each member of the pencil
of quintic threefolds
$$
Y_t:\qquad X_1X_2^4+X_2X_3^4+\ldots+X_5X_1^4\,-\,5tX_1X_2\cdots X_5\,=\,0
$$
is invariant under $H$. This leads, as above, to a splitting
$$
H^3(Y_t,\QQ)\,=\,T'_t\,\oplus\,S'_t,\qquad 
T'_t\,:=\,H^3(Y_t,\QQ)^H\,\cong\,\QQ^4,\qquad H^{3,0}(Y_t)\subset T'_t\otimes_\QQ\CC.
$$
Moreover, the degree four Picard-Fuchs equation defined by the 
variation of the Hodge structures $T'_t$ is the same as the Picard-Fuchs equation obtained from the variation of the  $T_t\cong H^3(M_t,\QQ)$.
In \cite{DGJ} more such examples are given. A possible explanation
would be that the CY threefold $M_t$ is birationally isomorphic to a desingularisation of $Y_t/H$. 

This is indeed the case.  Using results of Shioda, 
in the recent paper \cite{BGK} it is shown that there is a commutative 
diagram, where the arrows are rational maps which are quotients by certain
finite groups on suitable Zariski open subsets:
$$
\begin{array}{lcccr}
&&\tilde{X}_{dI,t}&&\\
&\swarrow&&\searrow&\\
X_t&&&&Y_t,\\
&\searrow&&\swarrow&\\
&&M_t&&
\end{array}
\qquad
\tilde{X}_{dI,t}:\quad
X_1^d+\ldots+X_5^d\,-\,5t(X_1X_2\cdots X_5)^{d/5}\,=\,0,
$$
where $d=5^2\cdot 41=1025$. Using the full diagram, one can show that 
the map $Y_t\ra M_t$ has degree $41$ and factors over $Y_t/H$. 
Thus there is a birational isomorphism
between $Y_t/H$ and $M_t$.

More generally, one can replace
$Y_t$ by a family with an equation
$$
\sum_{j=1}^5\,\prod_{i=1}^5\,X_i^{a_{ij}}\,-\,5tX_1\cdots X_5
$$
for suitable $5\times 5$ matrices with integer coefficients $a_{ij}$.
This again can be generalised to any number of variables.
A further generalisation to weighted projective spaces is given in \cite{B}.

\section{CY threefolds parametrised by Shimura varieties}

\subsection{} 
The period space $\cD=\cD_N$ parametrises the polarized weight three Hodge structures of CY-type on $(V_\ZZ\cong\ZZ^N,Q)$. Given a CY threefold $X$ and a marking, 
the unobstructedness of the deformations of $X$ implies that 
the period points of all deformations of $X$ are the points of a $q=(N-2)/2$-dimensional subvariety $B$ of $\cD$. 

On the other hand, there are many Hermitian symmetric domains which parametrise 
Hodge structures of CY-type. Such a domain is of the form $G(\RR)/K$,
where $G(\RR)$ is a real reductive Lie group which is the group of real points of an algebraic group defined over $\QQ$, and $K$ is maximal compact subgroup of $G(\RR)$.
Given $(V_\ZZ,Q)$,  there is a fixed representation of $G(\RR)$ on $V_\RR$
such that the image of $G(\RR)$ is contained in the symplectic group $Sp(Q,\RR)$. 
One considers the homomorphisms of real Lie groups
$$
h:\,S^1\,=\,\{z\in\CC:\,|z|=1\,\}\,\longrightarrow\,G(\RR)
$$
such that the eigenvalues of $h(z)$ are $z^p\bar{z}^q$ with
non-negative integers $p,q$ such that $p+q=3$.
Each such homomorphism gives a Hodge structure of weight three on $V_\ZZ$
by defining $V^{p,q}$ to be the eigenspace with eigenvalue $z^p\bar{z}^q$. 
The group $G(\RR)$ acts on the set of such Hodge structures by conjugation 
$h\mapsto ghg^{-1}$ with $g\in G(\RR)$. See for example \cite{Rlnm}, Chapter 1.

Changing the marking corresponds to an action of an element of $\Gamma:=Sp(V_\ZZ,Q)$ on $\cD$.
The moduli space of the Hodge structures of deformation of $X$ is thus the quotient of $B$ by the subgroup $\Gamma_B$ of $\Gamma$ which maps $B$ into itself. 
Not much is known about the subgroups $\Gamma_B$, 
see however \cite{DM}, \cite{R3} for a study in the case $N=4$.
In case $B$ is a Hermitian symmetric domain and $\Gamma_B$ is an arithmetic subgroup of $G(\QQ)$ one obtains a Shimura variety $\Gamma_B\backslash B=\Gamma_B\backslash G(\RR)/K$ which parametrises the deformations of $X$. 

Below we will review various examples from \cite{Rlnm}. 
In general it is not easy to decide if the deformations of a CY threefold are parametrised by a Shimura variety. See \cite{GMZ} for a family of CY threefolds
which are not parametrised by a Shimura variety.

\subsection{Example}(\cite{Bcycm}, $\S3$)\label{ex3elliptic-curve}
Let $E_i$, $i=1,2,3$ be elliptic curves and let $\iota_i:E_i\rightarrow E_i$
be the inversion $z\mapsto -z$ for the group law on $E_i$. Let 
$$
G_4\,:=\,\langle\,\iota_1\times\iota_2\times 1_{E_3},\;
\iota_1\times 1_{E_2} \times \iota_3\,\rangle\qquad\subset Aut(E_1\times E_2\times E_3).
$$
Then the (singular) variety $(E_1\times E_2\times E_3)/G_4$ has a resolution of singularieties which is a CY threefold $X$ with
$h^{2,1}=3$ (and $h^{1,1}=51$). Thus the deformation space of $X$ is
three dimensional. Obviously, it contains the CY varieties obtained by
deforming the three elliptic curves. Thus the period points of deformations of $X$ are in $B=\HH_1^3$, where $\HH_1$ is the upper half plane which parametrises elliptic curves.
Thus these CY's are parametrised by a Shimura variety.

\subsection{Examples of Borcea-Voisin type} 

\label{exBorceaVoisin}
Let $S$ be a K3 surface admitting an involution $\alpha_S$ such that 
$H^{2,0}(S)$ is in the eigenspace of the eigenvalue $-1$ for the action of $\alpha_S^*$ on $H^2(S,\CC)$.  We will assume moreover that the fixed locus of the involution $\alpha_S$ is made up of $k$ rational curves.
The dimension of the family of K3 surfaces admitting an involution acting non trivially on $H^{2,0}$ and fixing $k$ rational curves is 
$10-k$ and such a family is parametrised by a Shimura variety associated to
$SO(2,10-k)$.

Let $E$ be an elliptic curve and let $\iota$ be the involution 
$z\mapsto -z$ on $E$. 
The quotient threefold $(S\times E)/(\alpha_S\times \iota)$ admits a desingularisation which is a CY threefold $X$ 
(this construction is called Borcea-Voisin construction). 
In \cite{voisin}, \cite{B2} the Hodge numbers of $X$ are computed:
$$
h^{1,1}(X)\,=\,15+5k,\qquad h^{2,1}(X)\,=\,11-k.
$$
Hence the dimension of the family of the Calabi-Yau threefolds 
determined by $X$ is the sum of the dimension of the family of the K3 surfaces 
with involution
and the dimension of the family of elliptic curves.
Thus these CY threefolds are parametrised by the product of the 
Shimura varieties parametrising these two families, see \cite{Rlnm}, section 11.3.

The Example \ref{ex3elliptic-curve} is a particular case of this construction, indeed the desingularisation of the quotient 
$(E_1\times E_2)/(\iota_1\times\iota_2)$ is a K3 surface $S$ 
(in fact, it is a Kummer surface). The automorphism $\alpha_S$ induced on $S$ by $1_E\times \iota$ acts non trivially on $H^{2,0}(S)$ and fixes 
$8$ rational curves. 
Hence $(S\times E_3)/(\alpha_S\times \iota)$ is birational to 
$(E_1\times E_2\times E_3)/G_4$. 
For the Shimura varieties, one should remember that the real
Lie groups $SO(2,2)^0$ and $SL(2,\RR)\times SL(2,\RR)$ are isogeneous
and thus the Shimura variety associated to $SO(2,2)$ is indeed a quotient of 
$\HH_1\times\HH_1$.

\subsection{The easiest case}\label{easy}
Another particular case of the Borcea-Voisin construction is obtained 
by choosing $S$ to be the unique K3 surface with an automorphism $\alpha_S$ which fixes $k=10$ rational curves. 
This K3 surface $S$ is well known. It is described, for example, 
in \cite{SI} as the desingularisation of the quotient  
$$
(E_{\sqrt{-1}}\times E_{\sqrt{-1}})/(\gamma_{E}\times\gamma_{E}^3), 
\qquad\mbox{with}\quad
E_{\sqrt{-1}}\,:=\,\CC/(\ZZ+\sqrt{-1}\ZZ)
$$
and $\gamma_{E}$ is the automorphism of $E_{\sqrt{-1}}$ defined by
$z\mapsto \sqrt{-1}z$.

The third cohomology group of $X$ is:
$$
H^3(X,\QQ)\,\cong \, 
\left(T_S\otimes H^1(E,\QQ)\right)^{\alpha_S\times \iota}\,\cong\,
T_S\otimes H^1(E,\QQ),\qquad
T_S\,:=\,H^2(S,\QQ)^{\alpha_S=-1}\,\cong\,\QQ^2.
$$
As $S$, and thus the Hodge structure on $T_S$, is fixed, the variation of Hodge structures in $H^3(X,\ZZ)$ comes from the variation of Hodge structures
on $H^1(E,\ZZ)$. This rank $N=4$ variation is the direct sum of two (identical) rank two deformations, and is parametrised by $\HH_1$. In particular,
this variation of Hodge structures does not have maximally unipotent monodromy,
instead the monodromy operators have $2\times 2$ diagonal blocks in a suitable basis, cf.\ \cite{Rohde-maximal}, Example 2.7. We will see more examples
of variations without maximally unipotent monodromy in section \ref{mon}.

Any CY threefold from this family is also birationally isomorphic 
to a double cover of $\PP^3$ branched along the union of eight planes
As such it appears as entry no.\ 13 in the table in section 4.2.5 of the book \cite{M}. 
To obtain this double cover, one uses that $S$ is a double cover of $\PP^2$ branched over six lines. Putting one line `at infinity' in $\PP^2$,
a birational model of $S$ is (cf.\ \cite{GT}, 5.1, 5.2):
$$
S:\qquad s^2\,=\,xy(x-1)(y-1)(x-y).
$$
An elliptic curve $E$ can be defined by $t^2=u(u-1)(u-\lambda)$ for a
suitable $\lambda\in\CC$. Hence $X$, a desingularisation of $S\times E$ by the involution which fixes $x,y,u$ and maps $s,t\mapsto -s,-t$, has a birational model defined by
$$
X:\qquad w^2\,=\,xyu(x-1)(y-1)(u-1)(x-y)(u-\lambda).
$$
A suitable coordinate transformation on $\PP^3$ will map the banch locus to the one in Meyer's book \cite{M}.

\subsection{CY-type Hodge structures parametrised by Shimura varieties} 
There are many Shimura varieties which do parametrise variations of 
Hodge structures of CY-type, 
but where it is not known if these Hodge structures
come from CY threefolds.

For example, let $A$ be an abelian threefold and let $L\in H^2(A,\ZZ)$ be an ample divisor class. Note that $A$ is not a CY since $h^{1,0}=h^{2,0}=3$. One has $H^3(A,\ZZ)=\wedge^3H^1(A,\ZZ)$ 
and the primitive cohomology
$$
H^3(A,\ZZ)_{prim}\,\cong\,
H^3(A,\ZZ)/\left(L\wedge H^1(A,\ZZ)\right)\qquad(\cong\ZZ^{14})
$$ 
is a polarized
Hodge structure of $CY$ type with $q=h_{prim}^{2,1}=9-3=6$.
The moduli space $\cA_{3,L}$ of polarized abelian threefolds with the same polarization type as $(A,L)$ (these are all deformations of $A$) is the quotient  $\Gamma_L\backslash \HH_3$ of the Siegel half space by a discrete subgroup $\Gamma_L\subset Sp(6,\QQ)$. The image of $\HH_3$ in $\cD_{14}$ 
parametrizes the polarized Hodge structures of CY-type 
$H^3(A_t,\ZZ)_{prim}$ for deformations $(A_t,L_t)$ of $(A,L)$,
so these Hodge structures are parametrised by a Shimura variety.
To the best of our knowledge, it is not known if there exists a family of
CY threefolds $X_t$ such that
$H^3(A_t,\ZZ)_{prim}\cong H^3(X_t,\ZZ)/\mbox{torsion}$ for $t$ in an 
open, dense, subset of $\HH_3$.

A family of polarised CY-type Hodge structures, with $h^{2,1}=27$, parametrised by the Hermitian symmetric domain associated to the Lie group of type $E_7$ is
defined in \cite{Gross}. 
It is not yet known if there is a family of CY threefolds
with these Hodge structures.

\section{Examples with automorphisms of order three}\label{e3}

\subsection{Rohde's construction}\label{rohde}
In the paper \cite{Rohde-maximal}, J.C.\ Rohde constructs families of CY threefolds with $q=h^{2,1}=6-k$, for $0\leq k\leq 6$, which are parametrised by a $q$-dimensional Shimura variety, in this case a ball quotient. 
They are obtained as the desingularisation of the quotient of a product $E\times S$
by an automorphism of order three, where $E$ is a certain elliptic curve and $S$ is a K3 surface which admits an automorphism of order three which fixes $k$ rational curves and $k+3$ isolated points. These K3 surfaces were classified in \cite{AS}. A similar construction  with an automorphism of order four is discussed in \cite{Ga}, and various examples are given.

Let $\xi\in\CC$ be a primitive cube root of unity and consider the elliptic curve
$$
E:=\CC/\ZZ+\ZZ\xi,\qquad End(E)\,=\,\ZZ[\xi].
$$ 
We let $\alpha_E$ be the automorphism of $E$ defined by $z\mapsto \xi z$. 
A Weierstrass equation of $E$ is $y^2=x^3-1$ and 
$\alpha_E:(x,y)\mapsto (\xi x,y)$. 
The automorphism $\alpha_E$ gives the decomposition into eigenspaces, with eigenvalues $\xi$ and $\overline{\xi}$ respectively:
$$
H^1(E,\CC)\,=\,H^{1,0}(E)_\xi\,\oplus\,H^{0,1}(E)_{\overline{\xi}}.
$$

For any integer $k$ with $0\leq k\leq 6$, there exist K3 surfaces
$S$ with an automorphism of order three $\alpha_S$
such that the second cohomology group splits as:
$$
H^2(S,\QQ)\,=\,T_S\,\oplus_\perp\,N_S,\qquad 
N_S\,:=\,H^2(S,\QQ)^{\alpha_S}\cong \QQ^{8+2k},\quad
H^{2,0}(S)\subset T_S\otimes\CC.
$$
As $\dim H^2(S,\QQ)=22$, it follows that $\dim T_S=14-2k$
and the Hodge numbers of the weight two polarized Hodge structure
$T_S$ are $h^{2,0}(T_S)=1$, $h^{1,1}(T_S)=12-2k=2q$.
The action of $\alpha_S^*$ on $T_S$ defines a structure of $\QQ(\xi)$-vector space on  $T_S$. The eigenspaces for this action are

$$
T_S\otimes\CC\,=\,\underbrace{T^{2,0}_{S,\overline{\xi}}\,\oplus\, 
T^{1,1}_{S,\overline{\xi}}}_{T_{S,\overline{\xi}}}\,\oplus\,
\underbrace{T^{1,1}_{S,\xi}\,\oplus\,
T^{0,2}_{S,\xi}}_{T_{S,\xi}},\qquad 
\dim T^{1,1}_{S,\overline{\xi}}\,=\,
\dim T^{1,1}_{S,\xi}\,=\,6-k=q.
$$
The moduli space of such K3 surfaces is $q$-dimensional, and it is a quotient
of the $q$-ball in $\CC^q$, see section \ref{ball}, and it is a Shimura variety.

The weight three polarised rational Hodge substructure of 
$\alpha:=\alpha_S\times \alpha_E$-invariants in the tensor product
$H^2(S,\QQ)\otimes H^1(E,\QQ)$
is then of CY-type. Rohde shows that it is isomorphic to the third cohomology group 
of a CY threefold $X_S$ which is a desingularisation of the
(singular) quotient variety $(S\times E)/(\alpha_S\times \alpha_E)$:
$$
H^3(X_S,\QQ)\,\cong\,\left( H^2(S,\QQ)\otimes H^1(E,\QQ)\right)^\alpha
\,=\, \left(T_S\otimes  H^1(E,\QQ)\right)^{\alpha_S\times\alpha_E}.
$$
Here it is important that the fixed point locus of the automorphism $\alpha_S$ on $S$ consists of (smooth) rational curves and isolated points. 

The CY threefold $X_S$ still has an automorphism $\alpha_{X_S}$ 
of order three which is 
induced by $1_S\times \alpha_E$ (or, equivalently $\alpha_S^{-1}\times 1_E$).
As the eigenspaces of $\alpha_E$ on $H^1(E,\QQ)$ are 
$H^{1,0}(E)$ and $H^{0,1}(E)$ we obtain the decomposition into 
$\alpha_{X_S}$-eigenspaces:
$$
H^3(X_S,\CC)\,\cong\,\left( T_{S,\overline{\xi}}\otimes H^{1,0}(E)_\xi \right)
\,\oplus\, \left( T_{S,{\xi}}\otimes H^{0,1}(E)_{\overline{\xi}} \right)
\,=\,F^2\,\oplus\,\overline{F^2}
$$
where the last equality follows by inspection of the Hodge decomposition of $T_S$:
$$
H^{3,0}(X_S)\,\cong\,T^{2,0}_{S,\overline{\xi}}\otimes H^{1,0}(E)_\xi
\quad\subset\quad T_{S,\overline{\xi}}\otimes H^{1,0}(E)_\xi,
$$
and similarly
$$
H^{2,1}(X_S)\,\cong\,T^{1,1}_{S,\overline{\xi}}\otimes H^{1,0}(E)_\xi
\qquad\mbox{and thus}\quad
\dim H^{2,1}(X_S)\,=\,q.
$$
As the moduli of $S$ already provide a $q$-dimensional deformation 
space of $X_S$, one finds that all the deformations of $X_S$ 
are of this type.
Therefore these CY threefolds are parametrised by the same 
Shimura variety as the K3 surfaces $S$.

\subsection{No maximal unipotent monodromy}
\label{mon}
A peculiar feature of these families of CY threefolds is that
they do not have a large complex structure limit. In other 
words, their Picard-Fuchs equations do not have singular points with
maximally unipotent monodromy. To see why, recall that $F^2=H^3(X_S,\CC)_\xi$
and $\overline{F^2}=H^3(X_S,\CC)_{\overline{\xi}}$
are the eigenspaces of $\alpha_{X_S}$. The non-vanishing section $\omega$ of $F^3\subset F^2=F^2_\xi$ is always in the $\xi$-eigenspace,
which has dimension $1+q$. Therefore its derivatives under the Gauss Manin connection remain in this eigenspace. 
Instead of a degree $2(1+q)$ Picard-Fuchs equation one now finds an equation of degree $1+q$. 
This implies that the `standard' recipe for Mirror Symmetry cannot be applied to these families.  
In the case $q=1$, this equation was given explicitly in \cite{GG},
see also the next section. We already saw an another example of this in section \ref{easy}.

\subsection{The case $q=1$}\label{exR}
We recall our explicit description from \cite{GG} of the K3 surfaces from section
\ref{rohde} in the case the associated CY threefolds have 
$q=h^{2,1}=1$.

One starts with a polynomial $f=gh^2\in\CC[t]$ with $g,h$ of degree two
such that $f$ has four distinct zeroes, so up to the action of $Aut(\PP^1)$ we have only one parameter. The K3 surface $S_f$ has an elliptic fibration 
$\pi:S_f\rightarrow \PP^1$ with a section. 
Its Weierstrass model is:
$$
S_t:\quad Y^2=X^3+f(t)^2,\qquad f=gh^2,\quad\mbox{deg}(g)=\mbox{deg}(h)=2,
$$
where $t$ is the coordinate on $\PP^1$.
This surface has an automorphism of order three 
$$
\alpha_f:\;S_f\,\longrightarrow\,S_f,\qquad (X,Y,t)\,\longmapsto\,(\xi X,Y,t)
$$
which does act as $\overline\xi$  on $H^{2,0}(S_f)=\CC \rmd t\wedge\rmd X/Y$. 
Explicit computations show that 
$$
H^2(S_f,\QQ)\,=\,T_f\,\oplus\,N_f,\qquad 
N_f\,=\,H^2(S_f,\QQ)^{\alpha_f}\,\cong\,\QQ^{18},\quad 
H^{2,0}(S_f)\subset T_f\otimes\CC
$$
and that $\alpha_f$ fixes only five rational curves and eight isolated points.
The complexification of $T_f$ splits into four one-dimensional spaces:
$$
T_f\otimes\CC\,=\,T^{2,0}_{f,\overline{\xi}}\,\oplus\, T^{1,1}_{f,\overline{\xi}}\,\oplus\,
T^{1,1}_{f,{\xi}}\,\oplus\, T^{0,2}_{f,{\xi}}.
$$
Rohde's construction now produces a CY threefold $X_f$, the desingularisation of
$S_f\times E$ by the automorphism $\alpha=\alpha_f\times\alpha_E$ and
$$
H^3(X_f,\QQ)\,\cong\,(T_f\times H^1(E,\QQ))^\alpha \,=\,
\bigl(T^{2,0}_{f,\overline{\xi}}\,\oplus\, T^{1,1}_{f,\overline{\xi}}\bigr)\otimes H^{1,0}(E)_\xi\,\oplus\,\bigl(
T^{1,1}_{f,{\xi}}\,\oplus\, T^{0,2}_{f,{\xi}}\bigr)
\otimes H^{0,1}(E)_{\overline{\xi}}.
$$

The Hodge structures $T_f$ can be understood better by observing that all the smooth fibers of the elliptic fibration $\pi:S_f\rightarrow \PP^1$ are isomorphic (they are elliptic curves with $j$-invariant $0$, so are isomorphic to $E$).
Thus the elliptic fibration is isotrivial and becomes birationally isomorphic to a product 
after a base change. For this, we define a curve $C_f$ 
which is a 3:1 cyclic cover of $\PP^1$ with covering automorphism $\beta_f$:
$$
C_f:\quad v^3\,=\,f(t);\qquad \beta_f:\,C_f\rightarrow\,C_f,\quad 
(t,v)\,\longmapsto\,(t,\xi v).
$$
Substituting $f=v^3$ in the Weierstrass equation of $S_f$, one finds the birational isomorphism:
$$
C_f\times E\,\longrightarrow\,S_f\approx (C_f\times E)/(\beta_f\times \alpha_E),
\qquad ((t,v),(x,y))\,\longmapsto\,(X,Y,t)\,=\,(v^2x,v^3y,t).
$$
The automorphism $\alpha_f$ on $S_f$ is induced by $\alpha_E$.
This leads to an isomorphism of Hodge structures:
$$
T_f\,\cong\,\left(H^1(C_f,\QQ)\otimes H^1(E,\QQ)\right)^{\beta_f\times \alpha_E},
$$
which implies another isomorphism of Hodge structures:
$$
H^3(X_f,\QQ)\,\cong\,
\left(H^1(C_f,\QQ)\otimes H^1(E,\QQ)\otimes H^1(E,\QQ)\right)^H,
$$
where $H\cong (\ZZ/3\ZZ)^2$ is generated by the automorphisms 
$\beta_f\times \alpha_E\times 1_E$ and $1_{C_f}\times\alpha_E\times\alpha_E^{-1}$
of $C_f\times E\times E$.
This shows that the variation of the Hodge structures $H^3(X_f,\QQ)$ is entirely coming from the variation of the Hodge structures of the curves $C_f$.
Note that
$$
H^{3,0}(X_f)\,\cong \,H^{1,0}(C_f,\QQ)_{\overline{\xi}}\otimes 
H^1(E,\QQ)_\xi\otimes H^1(E,\QQ)_\xi.
$$

The Picard-Fuchs equations for the variation of Hodge structures of the curves $C_f$ is explicitly given in \cite{GG}. One can parametrise $\PP^1$ in such a way that 
$g(t)=t(t-1)$ and $h(t)=(t-\lambda)$ (and the other zero of $h$ is at $\infty$), thus
$C_f\cong C_\lambda$ with defining equation $v^3=t(t-1)(t-\lambda)^2$. 
The holomorphic one forms on this curve are 
$\rmd t/v$ and $(t-\lambda)\rmd t/v^2$, 
note that they have distinct eigenvalues 
$\overline{\xi},\xi$ for the automorphism $\beta_f$.
The Picard-Fuchs equation for $\eta:=\rmd t/v\in H^1(C_f,\QQ)_{\overline{\xi}}$ 
turns out to be:
$$
\left(\lambda(1-\lambda)\frac{\partial^2}{\partial \lambda^2}+
(1-2\lambda)\frac{\partial}{\partial \lambda}
-\frac{2}{3}\right)\eta\,=\,0.
$$
This is also the Picard-Fuchs equation for the holomorphic three form on the 
corresponding family of CY threefolds.

Rohde computes the Hodge numbers of these CY threefolds $X_f$ and finds:
$$
\dim H^{1,1}(X_f)\,=\,73,\qquad \dim H^{2,1}(X_f)\,=\,1.
$$
Any CY threefold $Y$ from the Mirror family, if it exists, should thus have
$h^{1,1}(Y)=1$ and $h^{2,1}(Y)=73$. At least three families of CY threefolds with these Hodge numbers are known: the complete intersections of type 
$(3,3)$ in $\PP^5$, $(2,2,3)$ in $\PP^6$ and $(4,4)$ in the weighted projective space
$\PP^5(1,1,1,1,2,2)$. But in these cases the Mirror families are known and they have maximally unipotent monodromy (cf. \cite{CYY}), hence they cannot be the Mirrors of the family of the $X_f$.

\subsection{The case $q>1$} \label{rohde_gen}
In case $q\leq 5$, we again find that the K3 surface $S$ has an isotrivial fibration
with smooth fibers isomorphic to $E$, 
but we could not find such a fibration in case $q=6$.
The CY threefold $X_S$ is then again a desingularisation of a quotient of the product of a curve $C$ with two copies of the fixed elliptic curve $E$. 
The variation of Hodge structures of the $X_S$ is obtained from the deformations of $C$.

For $q\leq 3$, we consider the surface (cf.\cite{GG})
$$
S_f: y^2=x^3+f(t)^2,\ \ f=gh^2,\qquad deg (f)=6,
$$
such that $g$ and $h$ have no common zeros and no multiple zeros. 
The curve $C_f:\ v^3=f(t)$ has the automorphism $\beta_f:(t,v)\ra (t,\xi v)$. 
As in \ref{exR}, Rohde's CY threefold $X_f$ is the desingularisation of $\left(C_f\times E\times E\right)/H$, 
where $H=\langle\beta_f\times \alpha_E\times 1_E, 
1_{C_f}\times \alpha_E\times\alpha_E^{-1}\rangle$ (see \cite[Remark 1.3]{GG}). 
The Hodge numbers of $X_f$ and the genus $g(C_f)$ of $C_f$ are as follows:
{\renewcommand{\arraystretch}{1.2}
$$
\begin{array}{|c|c|c|c|c|c}\hline
\mbox{deg}(g)&	\mbox{deg}(h)&g(C_f)&q=h^{2,1}(X_f)&h^{1,1}(X_f)\\
\hline
6&0&4&3&51\\
4&1&3&2&62\\
2&2&2&1&73\\
0&3&1&0&84\\ \hline
\end{array}
$$
}
The last line corresponds to a rigid CY threefold $X_f$ where $C_f\simeq E$, and $X_f$ is the desingularisation of the quotient $E\times E\times E$ by the group $\langle\alpha_E^{-1}\times \alpha_E\times 1_E, 1_{E}\times \alpha_E\times\alpha_E^{-1}\rangle$.
In this case, the K3 surface $S_f$ is described in \cite{SI}.

\

In case $q=4$, we consider the curve 
$$
C_l:\ \ v^6=l(t)\qquad \mbox{deg}(l)=12
$$ 
such that $l(t)$ has 5 double zeros.  It admits the automorphism $\beta_l:(t,v)\mapsto (t,\xi v)$. 
The quotient $(C_l\times E)/(\beta_l\times \alpha_E)$ is the K3 surface $S_l$ 
which has an elliptic fibration with Weierstrass equation 
$Y^2=X^3+l(t)$, where $X:=v^2x$, $Y:=v^3y$.\\
The desingularisation of $S_l$ is a K3 surface admitting an automorphism 
$\alpha_l$ of order 3 induced by $\alpha_E$. The fixed locus of $\alpha_f$
consists of $2$ rational curves and $5$ points. 
Applying Rhode's construction to the K3 surface $S_l$ one obtains a CY threefold $X$ such that $h^{2,1}(X)=q=4$ and $h^{1,1}(X)=40$. 

\

In case $q=5$, one needs a K3 surface $S$ with an automorphism $\alpha_S$
of order 3 which  fixes one rational curve and $4$ points (cf. \cite{Rohde-maximal}). In \cite{AS} a projective model of such a surface is given: 
it is a (singular) complete intersection in $\mathbb{P}^4$ with equations
$$
\left\{\begin{array}{l}F_2(x_0,\ldots x_3)=0,\\
G_3(x_0,\ldots x_3)\,=\,x_4^3,\end{array}\right.
$$
where $F_2$ and $G_3$ are homogeneous polynomials of degree 2 and 3 respectively. Moreover, the curve $V(F_2)\cap V(G_3)$ has $4$ singular points of type $A_1$. 
The surface $S$ is clearly a triple cover of 
the quadric defined $F_2=0$ in $\PP^3$ branched over the curve which is the intersection of this quadric with the cubic surface defined by $G_3=0$ 
in $\PP^3$. 
The inverse image in $S$ of a line in a ruling of the quadric in $\PP^3$ 
is an elliptic curve with a covering automorphism of order three which fixes the ramification points. Hence such an elliptic curve is isomorphic to $E$.
Thus $S$ admits an isotrivial fibration (in general without section)
with general fiber isomorphic to $E$. In this case Rohde's CY threefold has $h^{1,1}(X)=29$.

\subsection{The complex ball}\label{ball}
We briefly recall why the CY-type Hodge structures $H^3(X_S,\ZZ)$ are parametrised by a complex $q$-ball. More generally, with the notation from section \ref{pdom},
consider polarised weight three Hodge structures
on $(V_\ZZ\cong\ZZ^{2(1+q)},Q)$ of CY-type which, 
moreover, admit an automorphism
of order three $\phi$:
$$
\phi:\;V_\ZZ\,\longrightarrow\,V_\ZZ,\qquad 
Q(\phi x,\phi y)\,=\,Q(x,y),\qquad
\phi_\CC(V^{p,q})\,=\,V^{p,q},
\quad \phi^3=1_{V_\ZZ}.
$$

Then we have a decomposition of $V_\CC$ into $\phi$-eigenspaces, and we assume, as in the examples above, that the eigenspace of $\phi$ with eigenvalue $\xi$ is exactly $F^2$, so $F^2=F^2_\xi$. Then the $V^{p,q}$ are also $\phi$-eigenspaces: 
$$
V_\CC\,=\,V_\xi\,\oplus\,V_{\overline{\xi}}\,=\,
V^{3,0}_\xi\,\oplus\,V^{2,1}_\xi\,\oplus\,
V^{2,1}_{\overline{\xi}}\,\oplus\,V^{0,3}_{\overline{\xi}}.
$$
In particular, the subspace $F_2$, being an eigenspace of the fixed automorphism
$\phi$ of $V_\ZZ$, is now fixed in $V_\CC$. It remains to find the moduli of $V^{3,0}$ inside $F^2=V^{3,0}\oplus V^{2,1}$.
Recall the Hermitian form $H$ on $V_\CC$ which is positive definite on $V^{3,0}$
and negative definite on $V^{2,1}$. These two subspaces are perpendicular for $H$.
Thus the unitary group of $H_{|F^2}$ is isomorphic to the group $U(1,q)$.
It is well-known that this group acts transitively on the orthogonal decompositions
$F^2=W\oplus W^\perp$ with $H_{|W}>0$ (and thus $H_{|W^\perp}<0$).
The stabiliser of a given decomposition is the subgroup  $U(1)\times U(q)$,
hence the moduli space of these decompositions is the Hermitian symmetric domain
$$
U(1,q)/(U(1)\times U(q))\,\cong\,\BB^q\;=\;\{w\in \CC^q:\;|\!|w|\!|<1\,\}.
$$
It is easy to check that these decompositions correspond to the Hodge structures
under consideration, hence the $q$-ball is also the moduli space of these CY-type Hodge structures.

More explicitly, the Hermitian form $H$ has signature $(1+,q-)$ on 
the complex subspace $F^2\cong \CC^{1+q}$. Thus $F^2$ has a basis on which we have $H(z,z)=|z_0|^2-\sum_{j=1}^n|z_j|^2$.
In a Hodge decomposition, we must have $V^{3,0}=\CC w'$ for a non-zero
$w'=(w_0,w_1,\ldots,w_q)\in F^2$ such that $H(w',w')>0$, that is,
$|w_0|^2>\sum_{j=1}^n|w_j|^2$. In particular, $w_0\neq 0$ and so we 
may assume that $w_0=1$. Then $w'$ is determined by the point 
$w:=(w_1,\ldots, w_q)\in\CC^q$ with $\sum_{j=1}^n|w_j|^2<1$, 
that is, a point of the $q$-ball. Conversely, given $w\in \BB^q$, let $w'=(1,w)$
and define $V^{3,0}=\CC w'$, $V^{2,1}=(V^{3,0})^\perp$, the orthogonal complement, w.r.t.\ $H$, in $F^2$ of $V^{3,0}$ and
define $V^{1,2},V^{0,3}$ using $V^{p,q}=\overline{V^{q,p}}$. One easily checks that this gives a polarised Hodge structure on $(V_\ZZ,Q)$ which admits the automorphism $\phi$.

As we observed before in sections \ref{exR}, \ref{rohde_gen}, the ball also parametrises families of curves, like the $C_f$, and K3 surfaces, like the $S_f$.
Equivalently, it also parametrises certain Hodge structures of weight one and two.
The relation between these Hodge structures is given by the "half twist" construction,
see \cite{vG}, \cite{DK}.


\end{document}